\documentclass[12pt,reqno]{article}

\usepackage[usenames]{xcolor}
\usepackage{amsmath}
\usepackage{amssymb}
\usepackage{graphicx}
\usepackage{multirow}
\usepackage{breakurl}
\usepackage{tikz}
\usepackage{textcomp}

\definecolor{webgreen}{rgb}{0,.5,0}
\definecolor{webbrown}{rgb}{.6,0,0}

\usepackage[colorlinks=true,
linkcolor=webgreen,
filecolor=webbrown,
citecolor=webgreen]{hyperref}

\let\oeis\seqnum

\let\set\mathbb
\def\checkmark{\tikz\fill[scale=0.4](0,.35) -- (.25,0) -- (.7,.7) -- (.25,.15) -- cycle;}

\begin{document}

\begin{center}
  \vskip 1cm{\LARGE\bf
    How does the Gerrymander Sequence Continue?
  }
\vskip 1cm
\large
Manuel Kauers\\
Institute for Algebra\\
Johannes Kepler University\\
Altenberger Stra\ss e 69, 4040 Linz, Austria\\
\href{mailto:manuel.kauers@jku.at}{\tt manuel.kauers@jku.at}

\medskip
Christoph Koutschan\\
Johann Radon Institute for Computational and Applied Mathematics\\
Austrian Academy of Sciences\\
Altenberger Stra\ss e 69, 4040 Linz, Austria\\
\href{mailto:christoph.koutschan@oeaw.ac.at}{\tt christoph.koutschan@oeaw.ac.at}

\medskip
George Spahn\\
Department of Mathematics\\
Rutgers University (New Brunswick)\\
110 Frelinghuysen Road, Piscataway, NJ 08854-8019, USA\\
\href{mailto:geosp98@gmail.com}{\tt geosp98@gmail.com}
\end{center}

\vskip .2in

 \begin{abstract}
   We compute a few additional terms of the gerrymander sequence (OEIS
   \oeis{A348456}) and provide guessed equations for the generating functions
   of some sequences in its context.
 \end{abstract}

 \section{Introduction}

 In a guest lecture on April~28, 2022, in Zeilberger's
 famous course on experimental mathematics at Rutgers~\cite{zeilberger22}, Sloane
 talked about some of his favorite entries of the OEIS. One of the entries he
 highlighted in his lecture was \oeis{A348456}. The $n$th term of this sequence is defined as
 the number of ways to dissect a $2n\times2n$ chessboard into two polyominoes,
 each of area~$2n^2$. Here is one of the solutions for $n=4$:
 
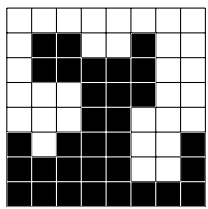
\begin{center}
 \begin{tikzpicture}[scale=.33]
  \draw(0,0) grid (8,8);
  \clip (0,0)--(8,0)--(8,3)--(7,3)--(7,1)--(5,1)--(5,4)--(6,4)--(6,7)--(5,7)--(5,6)--(3,6)--(3,7)--(1,7)--(1,5)--(3,5)--(3,3)--(2,3)--(2,2)--(1,2)--(1,3)--(0,3)--cycle;
  \fill[black](0,0) rectangle(8,8);
  \draw[white](0,0) grid(8,8);
 \end{tikzpicture}
 \end{center}
 Finding such dissections can be viewed as a combinatorial version of gerrymandering,
 and it has thus been suggested to call the sequence the gerrymander sequence.
 When Sloane gave his lecture, only the first three terms of the gerrymander sequence were known (they are $2$, $70$,~$80518$).
 He declared, perhaps exaggerating a bit, that he considers the next term of this sequence as
 the ``most wanted number'' in the whole OEIS.
 This statement motivated Zeilberger to offer a donation of \$100 to the OEIS in honor of the person who first manages
 to compute this most wanted number. In this paper, we explain how we computed not only the next term of \oeis{A348456}, but in fact
 the next four terms. They are
 \begin{alignat*}1
   &7157114189,\\
   &49852157614583644,\\
   &28289358593043414725944353,\text{ and}\\
   &1335056579423080371186456888543732162,
 \end{alignat*}
 respectively. In addition, we confirm the correctness of the previously known three terms.
 Our Mathematica source code is available at \url{http://www.koutschan.de/data/gerry/}.

 We employ the transfer-matrix method, a classical technique in enumerative combinatorics
 whose general idea is nicely explained in Sect.~4.7 of Stanley's textbook~\cite{stanley99}.
 We tweak this method by introducing catalytic variables, in order to apply it
 to the problem at hand; this is explained in Sect.~\ref{sec:transmat}.
 Besides the computation of specific terms, the transfer-matrix method also allows us to derive structural information
 about the generating function for boards of rectangular shapes $m\times n$ when $m$ is fixed and $n$ varies.
 The case $m=3$ was proposed by Knuth as a Monthly problem a few years ago~\cite{knuth18}.
 It turns out that if $a_n$ is the number of ways to break a $3\times2n$ board into two connected components of the same
 size (\oeis{A167242}), then
 \[
   \sum_{n=0}^\infty a_nx^n = \frac{1+\sqrt{1-4x}}{(\sqrt{1-4x}+x)^2}\frac1{\sqrt{1-4x}} - \frac{1-x^2+2x^3}{(1-x)^3}.
 \]
 It is not a coincidence but a consequence of a theorem of Furstenberg~\cite{furstenberg67} that the generating function is algebraic.
 In principle, it can be computed by combining the transfer-matrix method with the method of creative telescoping~\cite{Zeilberger91,Chyzak14,Koutschan13a}.
 However, this quickly becomes expensive when $m$ increases. In Sect.~\ref{sec:4}, we report on some computations we did
 in this direction.

 \section{The transfer matrix}\label{sec:transmat}

 The transfer-matrix method was invented in the
 context of statistical mechanics~\cite{KramersWannierI,KramersWannierII},
 in order to express the partition function
 of a statistical model in a simpler and more succinct form than its plain
 definition as a multi-dimensional sum.  The method is applicable whenever the
 mechanical system can be decomposed into a sequence of $N$ subsystems, each
 of them interacting only with the previous and the next one.  Let $\ell$
 denote the number of states that each of these subsystems can have, and
 $m_{i,j}(k)$ a ``statistical weight'' that is associated with state~$i$ of
 subsystem~$k-1$ being next to state~$j$ of subsystem~$k$. The relation
 between these two adjacent subsystems is then described by the \emph{transfer
   matrix} $M(k)=\bigl(m_{i,j}(k)\bigr)_{1\leq i,j\leq\ell}$, and the
 partition function of the whole system can be written in the form
 \[
   v_{\text{init}}^\top \cdot
   M(1) \cdot M(2) \cdots M(N) \cdot v_{\text{final}},
 \]
 where $v_{\text{init}}$ and $v_{\text{final}}$ are vectors of dimension~$\ell$.

 Recall that we are interested in counting the number of ways that an $m\times n$
 grid can be divided into two (or, more generally, $q$) connected regions, each of which
 is represented by a different color. We can apply the transfer-matrix method
 to this gerrymandering problem by decomposing the grid into a sequence of
 columns that are added one after the other. In each step, not all of the
 $2^m$ (resp., $q^m$) potential columns can be added, because some would
 violate the rules, e.g., by creating two disconnected regions of the same
 color:
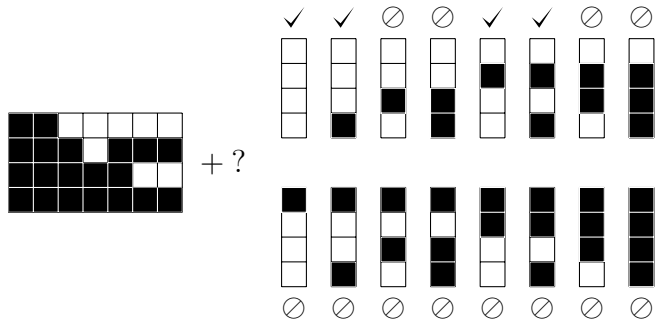
\begin{center}
\begin{tikzpicture}[scale=.33]
  \begin{scope}[xshift=-2cm]
   \draw(0,0) grid (7,4) (7,2) node[right] {${}+{?}$};
   \clip (0,0)--(7,0)--(7,1)--(5,1)--(5,2)--(7,2)--(7,3)--(4,3)--(4,2)--(3,2)--(3,3)--(2,3)--(2,4)--(0,4)--cycle;
   \fill[black] (0,0) rectangle(7,4);
   \draw[white](0,0) grid(7,4);
  \end{scope}
  \begin{scope}[yshift=3cm,xshift=9cm]
    \draw(0,0)grid(1,4) (.5,4) node[above] {\checkmark};
  \end{scope}
  \begin{scope}[yshift=3cm,xshift=11cm]
    \draw(0,0)grid(1,4) (.5,4) node[above] {\checkmark} ;
    \clip (0,0)rectangle(1,1);
    \fill[black](0,0) rectangle(1,4);
    \draw[white]grid(1,4);
  \end{scope}
  \begin{scope}[yshift=3cm,xshift=13cm]
    \draw(0,0)grid(1,4) (.5,4) node[above] {$\oslash$};
    \clip (0,1)rectangle(1,2);
    \fill[black](0,0) rectangle(1,4);
    \draw[white]grid(1,4);
  \end{scope}
  \begin{scope}[yshift=3cm,xshift=15cm]
    \draw(0,0)grid(1,4) (.5,4) node[above] {$\oslash$};
    \clip (0,0)rectangle(1,2);
    \fill[black](0,0) rectangle(1,4);
    \draw[white]grid(1,4);
  \end{scope}
  \begin{scope}[yshift=3cm,xshift=17cm]
    \draw(0,0)grid(1,4) (.5,4) node[above] {\checkmark};
    \clip (0,2)rectangle(1,3);
    \fill[black](0,0) rectangle(1,4);
    \draw[white]grid(1,4);
  \end{scope}
  \begin{scope}[yshift=3cm,xshift=19cm]
    \draw(0,0)grid(1,4) (.5,4) node[above] {\checkmark};
    \clip (0,0)rectangle(1,1) (0,2)rectangle(1,3);
    \fill[black](0,0) rectangle(1,4);
    \draw[white]grid(1,4);
  \end{scope}
  \begin{scope}[yshift=3cm,xshift=21cm]
    \draw(0,0)grid(1,4) (.5,4) node[above] {$\oslash$};
    \clip (0,1)rectangle(1,3);
    \fill[black](0,0) rectangle(1,4);
    \draw[white]grid(1,4);
  \end{scope}
  \begin{scope}[yshift=3cm,xshift=23cm]
    \draw(0,0)grid(1,4) (.5,4) node[above] {$\oslash$};
    \clip (0,0)rectangle(1,3);
    \fill[black](0,0) rectangle(1,4);
    \draw[white]grid(1,4);
  \end{scope}

  \begin{scope}[yshift=-3cm,xshift=9cm]
    \draw(0,0)grid(1,4) (.5,0) node[below] {$\oslash$};
    \clip (0,3)rectangle(1,4);
    \fill[black](0,0) rectangle(1,4);
    \draw[white]grid(1,4);
  \end{scope}
  \begin{scope}[yshift=-3cm,xshift=11cm]
    \draw(0,0)grid(1,4) (.5,0) node[below] {$\oslash$};
    \clip (0,0)rectangle(1,1) (0,3)rectangle(1,4);
    \fill[black](0,0) rectangle(1,4);
    \draw[white]grid(1,4);
  \end{scope}
  \begin{scope}[yshift=-3cm,xshift=13cm]
    \draw(0,0)grid(1,4) (.5,0) node[below] {$\oslash$};
    \clip (0,1)rectangle(1,2) (0,3)rectangle(1,4);
    \fill[black](0,0) rectangle(1,4);
    \draw[white]grid(1,4);
  \end{scope}
  \begin{scope}[yshift=-3cm,xshift=15cm]
    \draw(0,0)grid(1,4) (.5,0) node[below] {$\oslash$};
    \clip (0,0)rectangle(1,2) (0,3)rectangle(1,4);
    \fill[black](0,0) rectangle(1,4);
    \draw[white]grid(1,4);
  \end{scope}
  \begin{scope}[yshift=-3cm,xshift=17cm]
    \draw(0,0)grid(1,4) (.5,0) node[below] {$\oslash$};
    \clip (0,2)rectangle(1,4) ;
    \fill[black](0,0) rectangle(1,4);
    \draw[white]grid(1,4);
  \end{scope}
  \begin{scope}[yshift=-3cm,xshift=19cm]
    \draw(0,0)grid(1,4) (.5,0) node[below] {$\oslash$};
    \clip (0,0)rectangle(1,1) (0,2)rectangle(1,4);
    \fill[black](0,0) rectangle(1,4);
    \draw[white]grid(1,4);
  \end{scope}
  \begin{scope}[yshift=-3cm,xshift=21cm]
    \draw(0,0)grid(1,4) (.5,0) node[below] {$\oslash$};
    \clip (0,1)rectangle(1,4);
    \fill[black](0,0) rectangle(1,4);
    \draw[white]grid(1,4);
  \end{scope}
  \begin{scope}[yshift=-3cm,xshift=23cm]
    \draw(0,0)grid(1,4) (.5,0) node[below] {$\oslash$};
    \clip (0,0)rectangle(1,4);
    \fill[black](0,0) rectangle(1,4);
    \draw[white]grid(1,4);
  \end{scope}
 \end{tikzpicture}
 \end{center}

 Clearly, the information how the squares in the right-most column are colored
 is not sufficient to decide which columns can come next.  For example, we
 need to know that the two black squares in the last column belong to the same
 connected region, otherwise we could not add a column with only white
 squares.

 Thus, in order to decide which columns can be added and to tell whether the
 completed grid has the desired number of connected regions, we introduce
 \emph{states} for remembering connectivity information that is implied by
 all previous columns.  More precisely, a
 state is described by a pair $(c,P)$, where $c\in\{0,1\}^{2n}$ encodes the
 content of the last column (the colors white and black are now
 represented by the numbers $0$ and~$1$, respectively), and where
 $P=\{P_1,\dots,P_k\}$ is a partition of $\{1,\dots,2n\}$ that indicates
 which squares in that column belong to the same connected region. For
 example, in the figure above, the previous columns imply that the two black
 squares are connected, while the two white squares are not. The fact that
 they could (and should!) be connected by adding further columns is not
 relevant at this moment. Hence, in the partially completed grid we have
 $k$ different regions, and $P_j$ gives the positions of squares belonging to
 the $j$-th region. In particular, all these squares must have the same color,
 that is $|\{c_i\mid i\in P_j\}|=1$ for all~$1\leq j\leq k$.

 For example, $\bigl((0,0,1,0),\{\{1,2\},\{3,4\}\}\bigr)$ or
 $\bigl((1,1,0,0),\{\{1,2\},\{3\},\{4\}\}\bigr)$ are not valid state
 descriptions because the first violates the same-color condition, while the
 latter claims that the squares at positions $3$ and $4$ belong to different
 regions although they are obviously connected. In contrast,
 \begin{alignat*}1
   &\bigl((1,0,0,0),\{\{1\},\{2,3,4\}\}\bigr),\\
   &\bigl((1,0,0,1),\{\{1\},\{2,3\},\{4\}\}\bigr),\\
   \text{or }&\bigl((1,0,0,1),\{\{1,4\},\{2,3\}\}\bigr)
 \end{alignat*}
 are valid state descriptions, which we depict graphically as follows (connections
 that happen in previous columns are symbolized by an arc):
 
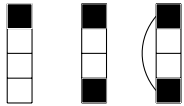
\begin{center}
  \begin{tikzpicture}[scale=.33]
    \begin{scope}
    \draw(0,0)grid(1,4);
    \clip(0,3)rectangle(1,4);
    \fill[black](0,0) rectangle(1,4);
    \draw[white]grid(1,4);
    \end{scope}
    \begin{scope}[xshift=3cm]
    \draw(0,0)grid(1,4);
    \clip(0,0)rectangle(1,1) (0,3)rectangle(1,4);
    \fill[black](0,0) rectangle(1,4);
    \draw[white]grid(1,4);
    \end{scope}
    \begin{scope}[xshift=6cm]
    \draw(0,0)grid(1,4);
    \draw(0,0.5)edge[bend left=40pt](0,3.5);
    \clip(0,0)rectangle(1,1) (0,3)rectangle(1,4);
    \fill[black](0,0) rectangle(1,4);
    \draw[white]grid(1,4);
    \end{scope}
  \end{tikzpicture}
 \end{center}
 But even if we stick to the above rules, there are still many pairs $(c,P)$
 that describe \emph{impossible} states, for example, $\bigl((0,1,0,1),\{\{1,3\},\{2,4\}\}\bigr)$.
 Connecting $1$ with~$3$ and $2$ with~$4$ produces arcs that cross each other,
 meaning that this connectivity cannot be achieved by extending the
 column to the left. Similarly, by considering how a column could be
 extended to the right, we can discard \emph{uninteresting}
 states, i.e., states that could possibly be reached, but which represent ``hopeless''
 situations that will never allow us to complete the grid in a satisfactory manner.
 For example, the state $\bigl((0,1,0,1),\{\{1\},\{2\},\{3\},\{4\}\}\bigr)$ is uninteresting,
 because all four squares are declared to belong to different regions, and it is
 impossible to connect $1$ with~$3$ and $2$ with~$4$ by adding more columns to the right.
 Hence, for the column $c=(0,1,0,1)$ we consider only two out of four possible states:
 
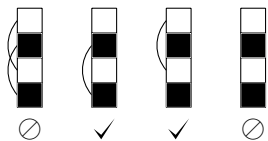
\begin{center}
  \begin{tikzpicture}[scale=.33]
    \begin{scope}
    \draw(0,0)grid(1,4) (.5,0) node[below] {$\oslash$};
    \draw(0,.5)edge[bend left=40pt](0,2.5);
    \draw(0,1.5)edge[bend left=40pt](0,3.5);
    \clip(0,0)rectangle(1,1) (0,2)rectangle(1,3);
    \fill[black](0,0) rectangle(1,4);
    \draw[white]grid(1,4);
    \end{scope}
    \begin{scope}[xshift=3cm]
    \draw(0,0)grid(1,4) (.5,0) node[below] {$\checkmark$};
    \draw(0,.5)edge[bend left=40pt](0,2.5);
    \clip(0,0)rectangle(1,1) (0,2)rectangle(1,3);
    \fill[black](0,0) rectangle(1,4);
    \draw[white]grid(1,4);
    \end{scope}
    \begin{scope}[xshift=6cm]
    \draw(0,0)grid(1,4) (.5,0) node[below] {$\checkmark$};
    \draw(0,1.5)edge[bend left=40pt](0,3.5);
    \clip(0,0)rectangle(1,1) (0,2)rectangle(1,3);
    \fill[black](0,0) rectangle(1,4);
    \draw[white]grid(1,4);
    \end{scope}
    \begin{scope}[xshift=9cm]
    \draw(0,0)grid(1,4) (.5,0) node[below] {$\oslash$};
    \clip(0,0)rectangle(1,1) (0,2)rectangle(1,3);
    \fill[black](0,0) rectangle(1,4);
    \draw[white]grid(1,4);
    \end{scope}
  \end{tikzpicture}
 \end{center}

 Algorithmically, we construct the set of states as follows: in step (1) we
 enumerate the set of all states that are not impossible, and in step (2) we
 discard those states which are uninteresting. In both steps it is clear from
 the construction that no necessary state is discarded, ensuring the
 correctness of our method.
 \begin{enumerate}
 \item For each tuple $c\in\{0,1\}^{2n}$, all non-crossing arc configurations
   are generated, where the vertices for these configurations are maximal
   chunks of squares of the same color. All arcs point in the same direction
   and each arc connects two vertices of the same color. Each vertex has at
   most one outgoing and at most one incoming arc. For example, the following
   column of size~8 admits seven such arc configurations:
   
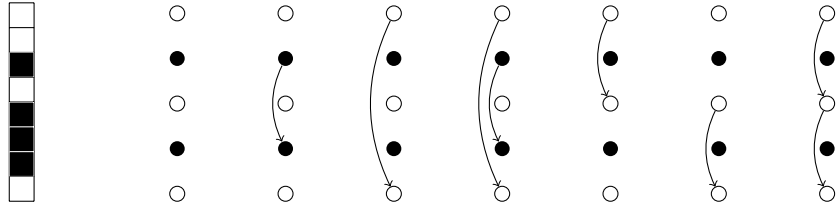
\begin{center}   
   \begin{tikzpicture}[scale=0.33]
     \draw(0,0)grid(1,8) (.5,0);
     \clip(0,1)rectangle(1,4) (0,5)rectangle(1,6);
     \fill[black](0,0) rectangle(1,8);
     \draw[white]grid(1,8);
   \end{tikzpicture}
   \qquad\qquad
   \begin{tikzpicture}[scale=0.6,xscale=.8]
     \foreach \x in {0,...,6} {
       \foreach \i in {1,3,5} \node[circle, inner sep=2pt, draw] (v\x\i) at (3*\x,5-\i) {};
       \foreach \i in {2,4} \node[circle, inner sep=2pt, fill] (v\x\i) at (3*\x,5-\i) {};
     }
     \draw[->] (v12) to[out=240,in=120] (v14);
     \draw[->] (v21) to[out=240,in=120] (v25);
     \draw[->] (v31) to[out=240,in=120] (v35);
     \draw[->] (v32) to[out=240,in=120] (v34);
     \draw[->] (v41) to[out=240,in=120] (v43);
     \draw[->] (v53) to[out=240,in=120] (v55);
     \draw[->] (v61) to[out=240,in=120] (v63);
     \draw[->] (v63) to[out=240,in=120] (v65);
   \end{tikzpicture}
   \end{center}
 \item A state $(c,P)$ is uninteresting if there are two connected components
   $A_1,A_2\in P$ of the same color and $B_1,B_2\in P$ of the opposite color
   such that the following condition holds:
   \begin{alignat*}1
     &\max A_1 + 1 = \min B_1 \land \max B_1 = \min A_2 - 1\\
     &\land(\max A_2 < \max B_2 \lor \min B_2 < \min A_1).
   \end{alignat*}
   It corresponds to situations where we cannot draw another non-crossing arc
   configuration on the right side such that all vertices of the same color get
   connected:
   
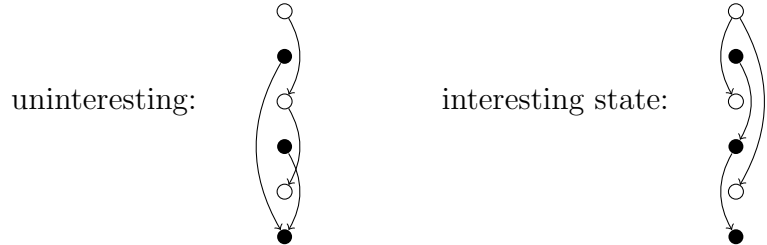
\begin{center}
   \begin{tikzpicture}[scale=0.6]
     \foreach \x in {1,2} {
       \foreach \i in {1,3,5} \node[circle, inner sep=2pt, draw] (v\x\i) at (10*\x,6-\i) {};
       \foreach \i in {2,4,6} \node[circle, inner sep=2pt, fill] (v\x\i) at (10*\x,6-\i) {};
     }
     \node at (6,3) {uninteresting:};
     \node at (16,3) {interesting state:};
     \draw[->] (v12) to[out=240,in=120] (v16);
     \draw[->] (v11) to[out=300,in=60] (v13);
     \draw[->] (v13) to[out=300,in=60] (v15);
     \draw[->] (v14) to[out=300,in=60] (v16);
     \draw[->] (v21) to[out=240,in=120] (v23);
     \draw[->] (v24) to[out=240,in=120] (v26);
     \draw[->] (v22) to[out=300,in=60] (v24);
     \draw[->] (v21) to[out=300,in=60] (v25);
   \end{tikzpicture}
   \end{center}
 \end{enumerate}
 
 There is one subtle issue one still has to take care of: if the current
 column equals $(0,\dots,0)$ or $(1,\dots,1)$, then we need to store the
 information whether the other color has not yet appeared (in which case this
 column can be followed by any other column), or whether the other color has
 appeared previously (in which case this column can only be followed by more
 copies of the same column).  In our graphical notation, we decorate the
 latter of these two states by a prime.

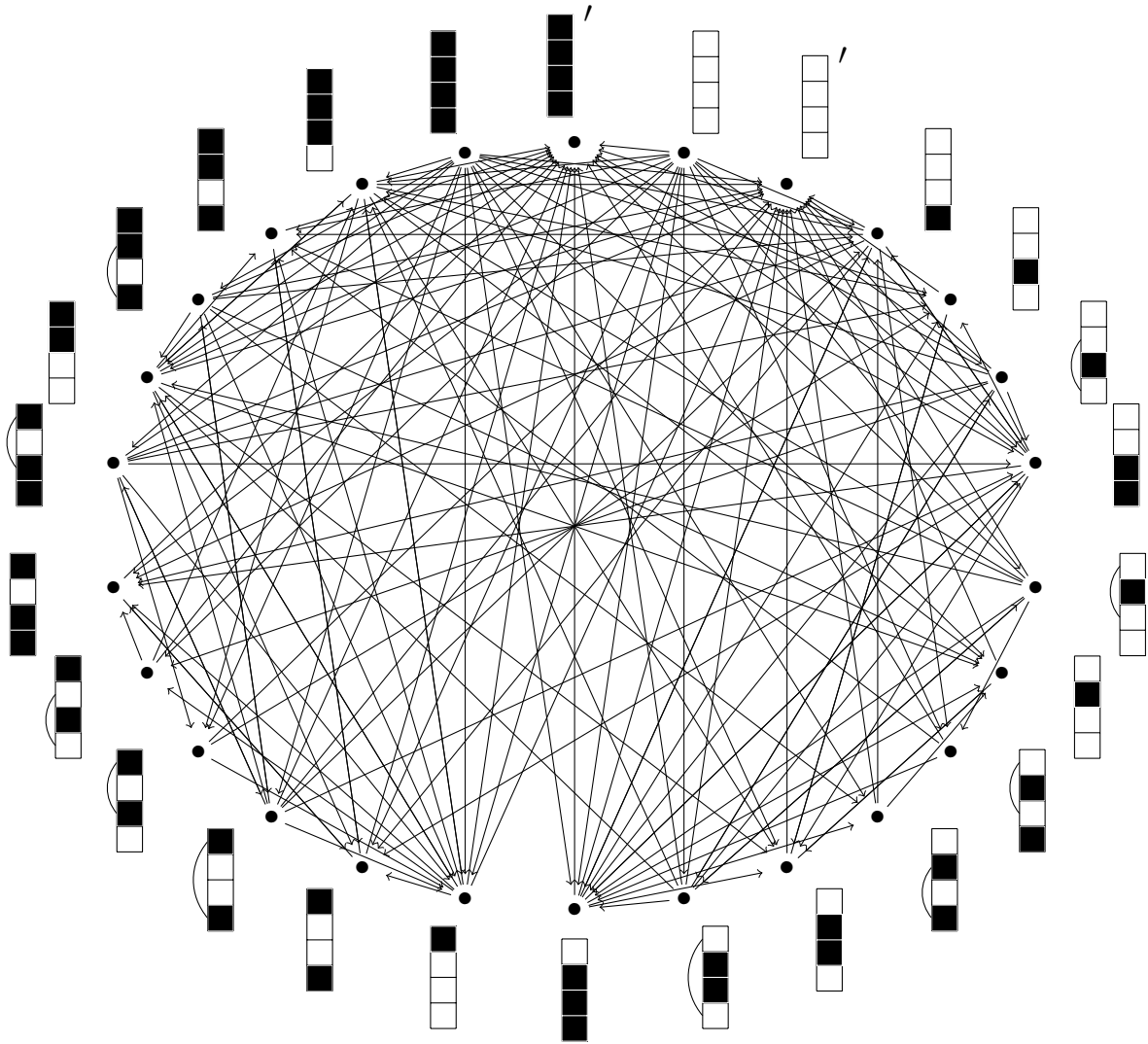
\begin{figure}
 \begin{center}
  \def\a{\begin{tikzpicture}[scale=.35]\draw(0,0)grid(1,4);}\def\b{\fill[black](0,0)rectangle(1,4);\draw[white](0,0)grid(1,4);\end{tikzpicture}}
  \begin{tikzpicture}[scale=5, xscale=1.25]
    \begin{scope}[scale=1.2]
    \node (1a) at (0.239316, 0.970942) {\a\clip(0,0);\b};
    \node (2a) at (0.464723, 0.935456) {\a\draw(1,4)node[right]{$\prime$};\clip(0,0);\b};
    \node (3a) at (0.663123, 0.748511) {\a\clip(0,0)rectangle(1,1);\b};
    \node (4a) at (0.822984, 0.568065) {\a\clip(0,1)rectangle(1,2);\b};
    \node (5a) at (0.935016, 0.354605) {\a\draw(0,.5)edge[bend left=40pt](0,2.5);\clip(0,1)rectangle(1,2);\b};
    \node[xshift=.1cm] (6a) at (0.992709, 0.120537) {\a\clip(0,0)rectangle(1,2);\b};
    \begin{scope}[yshift=-.1cm]    
    \node[xshift=.1cm] (7a) at (0.992709, -0.120537) {\a\draw(0,1.5)edge[bend left=40pt](0,3.5);\clip(0,2)rectangle(1,3);\b};
    \node (8a) at (0.935016, -0.354605) {\a\clip(0,2)rectangle(1,3);\b};
    \node (9a) at (0.822984, -0.568065) {\a\draw(0,1.5)edge[bend left=40pt](0,3.5);\clip(0,0)rectangle(1,1)(0,2)rectangle(1,3);\b};
    \node (10a) at (0.663123, -0.748511) {\a\draw(0,.5)edge[bend left=40pt](0,2.5);\clip(0,0)rectangle(1,1)(0,2)rectangle(1,3);\b};
    \node (11a) at (0.464723, -0.885456) {\a\clip(0,1)rectangle(1,3);\b};
    \node (12a) at (0.239316, -0.970942) {\a\draw(0,.5)edge[bend left=40pt](0,3.5);\clip(0,1)rectangle(1,3);\b};
    \node (13a) at (0., -1.) {\a\clip(0,0)rectangle(1,3);\b};
    \node (14a) at (-0.239316, -0.970942) {\a\clip(0,3)rectangle(1,4);\b};
    \node (15a) at (-0.464723, -0.885456) {\a\clip(0,0)rectangle(1,1)(0,3)rectangle(1,4);\b};
    \node (16a) at (-0.663123, -0.748511) {\a\draw(0,.5)edge[bend left=40pt](0,3.5);\clip(0,0)rectangle(1,1)(0,3)rectangle(1,4);\b};
    \node (17a) at (-0.822984, -0.568065) {\a\draw(0,1.5)edge[bend left=40pt](0,3.5);\clip(0,1)rectangle(1,2)(0,3)rectangle(1,4);\b};
    \node (18a) at (-0.935016, -0.354605) {\a\draw(0,.5)edge[bend left=40pt](0,2.5);\clip(0,1)rectangle(1,2)(0,3)rectangle(1,4);\b};
    \node[xshift=-.1cm] (19a) at (-0.992709, -0.120537) {\a\clip(0,0)rectangle(1,2)(0,3)rectangle(1,4);\b};
    \end{scope}
    \node[xshift=-.1cm] (20a) at (-0.992709, 0.120537) {\a\draw(0,1.5)edge[bend left=40pt](0,3.5);\clip(0,0)rectangle(1,2)(0,3)rectangle(1,4);\b};
    \node (21a) at (-0.935016, 0.354605) {\a\clip(0,2)rectangle(1,4);\b};
    \node (22a) at (-0.822984, 0.568065) {\a\draw(0,.5)edge[bend left=40pt](0,2.5);\clip(0,0)rectangle(1,1)(0,2)rectangle(1,4);\b};
    \node (23a) at (-0.663123, 0.748511) {\a\clip(0,0)rectangle(1,1)(0,2)rectangle(1,4);\b};
    \node (24a) at (-0.464723, 0.885456) {\a\clip(0,1)rectangle(1,4);\b};
    \node (25a) at (-0.239316, 0.970942) {\a\clip(0,0)rectangle(1,4);\b};
    \node (26a) at (0., 1.03) {\a\draw(1,4)node[right]{$\prime$};\clip(0,0)rectangle(1,4);\b};
    \end{scope}

    \node[circle] (1) at (0.239316, 0.970942) {$\bullet$};
    \node[circle] (2) at (0.464723, 0.885456) {$\bullet$};
    \node[circle] (3) at (0.663123, 0.748511) {$\bullet$};
    \node[circle] (4) at (0.822984, 0.568065) {$\bullet$};
    \node[circle] (5) at (0.935016, 0.354605) {$\bullet$};
    \node[circle, xshift=.1cm] (6) at (0.992709, 0.120537) {$\bullet$};
    \begin{scope}[yshift=-.1cm]    
    \node[circle, xshift=.1cm] (7) at (0.992709, -0.120537) {$\bullet$};
    \node[circle] (8) at (0.935016, -0.354605) {$\bullet$};
    \node[circle] (9) at (0.822984, -0.568065) {$\bullet$};
    \node[circle] (10) at (0.663123, -0.748511) {$\bullet$};
    \node[circle] (11) at (0.464723, -0.885456) {$\bullet$};
    \node[circle] (12) at (0.239316, -0.970942) {$\bullet$};
    \node[circle] (13) at (0., -1.) {$\bullet$};
    \node[circle] (14) at (-0.239316, -0.970942) {$\bullet$};
    \node[circle] (15) at (-0.464723, -0.885456) {$\bullet$};
    \node[circle] (16) at (-0.663123, -0.748511) {$\bullet$};
    \node[circle] (17) at (-0.822984, -0.568065) {$\bullet$};
    \node[circle] (18) at (-0.935016, -0.354605) {$\bullet$};
    \node[circle, xshift=-.1cm] (19) at (-0.992709, -0.120537) {$\bullet$};
    \end{scope}
    \node[circle, xshift=-.1cm] (20) at (-0.992709, 0.120537) {$\bullet$};
    \node[circle] (21) at (-0.935016, 0.354605) {$\bullet$};
    \node[circle] (22) at (-0.822984, 0.568065) {$\bullet$};
    \node[circle] (23) at (-0.663123, 0.748511) {$\bullet$};
    \node[circle] (24) at (-0.464723, 0.885456) {$\bullet$};
    \node[circle] (25) at (-0.239316, 0.970942) {$\bullet$};
    \node[circle] (26) at (0., 1.) {$\bullet$};    

    \node[circle] (1o) at (0.239316, 0.970942) {};
    \node[circle] (2o) at (0.464723, 0.885456) {};
    \node[circle] (3o) at (0.663123, 0.748511) {};
    \node[circle] (4o) at (0.822984, 0.568065) {};
    \node[circle] (5o) at (0.935016, 0.354605) {};
    \node[circle, xshift=.1cm] (6o) at (0.992709, 0.120537) {};
    \begin{scope}[yshift=-.1cm]    
    \node[circle, xshift=.1cm] (7o) at (0.992709, -0.120537) {};
    \node[circle] (8o) at (0.935016, -0.354605) {};
    \node[circle] (9o) at (0.822984, -0.568065) {};
    \node[circle] (10o) at (0.663123, -0.748511) {};
    \node[circle] (11o) at (0.464723, -0.885456) {};
    \node[circle] (12o) at (0.239316, -0.970942) {};
    \node[circle] (13o) at (0., -1.) {};
    \node[circle] (14o) at (-0.239316, -0.970942) {};
    \node[circle] (15o) at (-0.464723, -0.885456) {};
    \node[circle] (16o) at (-0.663123, -0.748511) {};
    \node[circle] (17o) at (-0.822984, -0.568065) {};
    \node[circle] (18o) at (-0.935016, -0.354605) {};
    \node[circle, xshift=-.1cm] (19o) at (-0.992709, -0.120537) {};
    \end{scope}
    \node[circle, xshift=-.1cm] (20o) at (-0.992709, 0.120537) {};
    \node[circle] (21o) at (-0.935016, 0.354605) {};
    \node[circle] (22o) at (-0.822984, 0.568065) {};
    \node[circle] (23o) at (-0.663123, 0.748511) {};
    \node[circle] (24o) at (-0.464723, 0.885456) {};
    \node[circle] (25o) at (-0.239316, 0.970942) {};
    \node[circle] (26o) at (0., 1.) {};    
    
    \foreach\i in {3,5,6,7,9,12,13,14,15,18,19,21,23,24,26} \draw[->] (1o)edge(\i);
    \foreach\i in {2,6,9,13,15,19,23,26} \draw[->] (3o)edge(\i);
    \foreach\i in {2,11} \draw[->] (4o)edge(\i);
    \foreach\i in {2,6,12,13,18,19,24,26} \draw[->] (5o)edge(\i);
    \foreach\i in {2,3,4,11,13,15,19,26} \draw[->] (6o)edge(\i);
    \foreach\i in {2,9,12,13,21,23,24,26} \draw[->] (7o)edge(\i);
    \foreach\i in {2,11} \draw[->] (8o)edge(\i);
    \foreach\i in {13,23,26} \draw[->] (9o)edge(\i);
    \foreach\i in {2,3,8} \draw[->] (10o)edge(\i);
    \foreach\i in {2,4,8} \draw[->] (11o)edge(\i);
    \foreach\i in {2,5,6,7,13,21,24,26} \draw[->] (12o)edge(\i);
    \foreach\i in {2,3,4,6,8,10,11,26} \draw[->] (13o)edge(\i);
    \foreach\i in {2,15,18,19,21,23,24,26} \draw[->] (14o)edge(\i);
    \foreach\i in {19,23,26} \draw[->] (15o)edge(\i);
    \foreach\i in {2,3,6,14,20,21,22,26} \draw[->] (16o)edge(\i);
    \foreach\i in {2,4,14} \draw[->] (17o)edge(\i);
    \foreach\i in {19,24,26} \draw[->] (18o)edge(\i);
    \foreach\i in {15,26} \draw[->] (19o)edge(\i);
    \foreach\i in {2,3,4,6,14,16,17,26} \draw[->] (20o)edge(\i);
    \foreach\i in {2,8,11,14,15,23,24,26} \draw[->] (21o)edge(\i);
    \foreach\i in {2,3,8,10,14,16,21,26} \draw[->] (22o)edge(\i);
    \foreach\i in {15,26} \draw[->] (23o)edge(\i);
    \foreach\i in {2,4,8,11,14,17,21,26} \draw[->] (24o)edge(\i);
    \foreach\i in {2,3,4,6,8,10,11,13,14,16,17,20,21,22,24} \draw[->] (25o)edge(\i);
  \end{tikzpicture}
  \end{center}
  \caption{All 26 states for a grid with 4 rows and their possible transitions;
    note that each state can also be followed by itself --- these loops are not
    depicted.}
  \label{fig:transitions}
 \end{figure}

 Let us denote by~$L$ the set of states that is constructed according to the
 above rules. The number of states~$\ell=|L|$  grows (at least)
 exponentially with~$n$, since each of the $2^{2n}$ possible columns appears
 in at least one state. For example, for $n=2$ we have 16 different columns
 but 26 states in total; they are explicitly enumerated in
 Figure~\ref{fig:transitions}. The number of states for $1\leq n\leq7$ is
 given in Table~\ref{tab:sparsity}.

 To construct the transfer matrix~$M$, for each state we need to determine
 which other states it can move into by adding an appropriate next column. Note that
 for this purpose it does not matter in which particular column of the grid we
 are: in each step we can use the same matrix~$M$, which is in contrast to the
 general situation sketched at the beginning of this section. The rows and
 columns of the transfer matrix are indexed by the states; hence we obtain
 an $\ell\times\ell$ matrix. Let $s=(c,P)$ be any of the $\ell$ states and
 let $c'\in\{0,1\}^{2n}$ be an arbitrary column. If attaching $c'$ to the
 state~$s$ would violate the connectivity requirements, then the matrix
 entries at positions $(s,s')$ are set to~$0$, where $s'$ is any state
 containing~$c'$.

 But what should be put as the matrix entry when a transition is actually
 possible? Here another condition has to be considered that so far has not
 been taken care of: the two regions must have the same area. Since this
 requirement can only be checked at the very end when the whole grid is
 filled, we need to propagate information about the number of squares of
 either color through the whole computation. For this purpose, we introduce a
 ``catalytic'' variable~$x$ that counts the number of white squares that have
 been used so far.  In each transition this counter has to be increased
 accordingly.
 Assume that state~$s=(c,P)$ admits attaching column $c'$ to it, which
 basically means that no existing region gets disconnected (this happens when
 $c$ and $c'$ have opposite colors at all positions given by some $P_j\in P$),
 except when all squares of $c'$ have the same color and only a single $P_j\in P$
 is related to the other color, in which case we move to one of the
 two primed states. In any case, the new connectivity information $P'$ is
 uniquely determined from $c,c',P$, and therefore yields a new state
 $s'=(c',P')$. It may well happen that $s'\not\in L$ is an uninteresting state,
 in which case no matrix entry is generated. Otherwise the matrix entry at
 position $(s,s')$ is set to $x^{\#_0(c')}$, where $\#_0(c')$ denotes the
 number of $0$'s in~$c'$.

 \begin{table}
   \begin{center}
     \begin{tabular}{r|r|r|r@{\,\%\ }|r|r@{\,\%\ }|r|r@{\,\%\ }}
       \multirow{2}{*}{$n$} & \multirow{2}{*}{$\quad\ \ell\quad\ $} &
       \multicolumn{2}{c|}{$v_{\text{init}}$} &
       \multicolumn{2}{c|}{$v_{\text{final}}$} &
       \multicolumn{2}{c}{$M$} \\
       \cline{3-8} \rule{0pt}{12pt} & &
       $\qquad\quad\#_{\neq0}$ & \multicolumn{1}{r|}{sparsity} &
       $\qquad\quad\#_{\neq0}$ & \multicolumn{1}{r|}{sparsity} &
       $\qquad\quad\#_{\neq0}$ & \multicolumn{1}{r}{sparsity} \\
       \hline \rule{0pt}{12pt}
       1 & 6      & 4   & 33.33 & 6   &  0.00 & 16        & 55.56 \\
       2 & 26     & 14  & 46.15 & 16  & 38.46 & 178       & 73.67 \\
       3 & 154    & 32  & 79.22 & 34  & 77.92 & 2546      & 89.26 \\
       4 & 1026   & 58  & 94.35 & 60  & 94.15 & 44008     & 95.82 \\
       5 & 7222   & 92  & 98.73 & 94  & 98.70 & 832454    & 98.40 \\
       6 & 52650  & 134 & 99.75 & 136 & 99.74 & 16505486  & 99.40 \\
       7 & 393878 & 184 & 99.95 & 186 & 99.95 & 337332580 & 99.78
     \end{tabular}
   \end{center}
   \caption{Number $\ell$ of states and number $\#_{\neq0}$ of non-zero entries
     in $v_{\text{init}}$, $v_{\text{final}}$, and in the transfer matrix~$M$.}
   \label{tab:sparsity}
 \end{table}

 Now that we have constructed the transfer matrix~$M$, it remains to consider
 the start and the end of this process. We start the grid with a single
 column. There cannot be any additional connectivity information other than
 what can be seen in this column. We define a start vector $v_{\text{init}}$, which is
 indexed by the states and hence is $\ell$-dimensional. Its entry at
 position $s=(c,P)$ equals $x^{\#_0(c)}$ if the parts of~$P$ correspond
 exactly to the consecutive runs of entries of the same color in~$c$ (in other
 words, if the graphical representation of $s$ has no arcs (and no prime!)),
 and $0$ otherwise.

 When we have filled the grid up to the last column, we have to decide which
 states are ``accept'' states. Clearly, this is the case for states $s=(c,P)$
 such that $|P|\leq2$. Since we just want to add up the results of all
 acceptable states, we define a vector $v_{\text{final}}$ that is $1$ at accept states
 and $0$ otherwise (see Table~\ref{tab:sparsity}).

 Having all this at hand---the transfer matrix~$M$, the start vector~$v_{\text{init}}$,
 and the end vector~$v_{\text{final}}$---one can compute
 \begin{equation}\label{eq:vMv}
   p(x) = v_{\text{init}}^\top \cdot M^{2n-1} \cdot v_{\text{final}},
 \end{equation}
 which is a polynomial in~$x$. The coefficient of $x^k$ in $p$ gives the
 number of ways that the $2n\times2n$ grid can be divided into a white
 polyomino of $k$ squares and a black polyomino consisting of $4n^2-k$
 squares. This means that for sequence \oeis{A348456} we need to extract the
 coefficient of $x^{2n^2}$ and divide the result by $2$ in order to
 eliminate the distinction between black and white.  

 \section{Optimizations}\label{sec:3}

 In this section, we discuss some optimizations that one can employ when
 implementing the method described in Section~\ref{sec:transmat}. For our
 purposes, we have used the computer algebra systems Maple and Mathematica.

 First, it is clear that in~\eqref{eq:vMv} one should not compute $M^{2n-1}$
 explicitly, using expensive matrix multiplications. Instead, it is more
 efficient to exploit the associativity of matrix multiplication and perform
 only (cheap) matrix-vector multiplications:
 \[
   p(x) = (\cdots((v_{\text{init}}^\top\cdot M)\cdot M)\cdots)\cdot v_{\text{final}}.
 \]

 Then we address the problem that the matrices get large: for example, for
 $n=7$ the transfer matrix~$M$ has $393878^2$, i.e., about $115$ billion entries,
 which poses memory challenges when one does not have a supercomputer at hand.
 However, we realize that the matrix is very sparse. Clearly, in each
 row we can have at most $2^{2n}$ nonzero entries (because this is the number
 of different grid columns we can add), but actually there are much fewer,
 since many of these columns are not admissible (e.g., because they disconnect
 existing regions). Table~\ref{tab:sparsity} shows the sparsity of the
 transfer matrix~$M$ and of the vectors $v_{\text{init}}$ and $v_{\text{final}}$.
 In Mathematica, one can use the command \texttt{SparseArray} to store such
 matrices in a memory-efficient way, which has the additional advantage that
 it also speeds up the matrix-vector multiplications. Note that the vector
 will quickly become dense as we multiply the matrix to it, reflecting the fact
 that we have no unreachable states in~$L$.

 The next observation concerns the structure of the transfer matrix. Since
 each specific column of the matrix is responsible for producing the counting
 polynomial for a specific state $(c,P)$, by matrix multiplication, all
 entries of this column of~$M$ must either be~$0$ or~$x^{\#_0(c)}$.  No other
 powers of $x$ can occur in the same column. It is more efficient to work with
 the $\{0,1\}$-matrix $M'=M\big|_{x\to1}$ and store the $x$-powers in a
 separate diagonal matrix
 \[
   X = \operatorname{diag}\bigl(\bigl( x^{\#_0(c)} \bigr)_{(c,P)\in L}\bigr)
 \]
 such that $M=M'\cdot X$. Computing $v_{\text{init}}^\top\cdot M$ has the disadvantage that
 large intermediate expressions are produced that can only be combined after
 expansion, since they are sums of products of monomials times polynomials.
 This is avoided by computing $\bigl(v_{\text{init}}^\top\cdot M'\bigr)\cdot X$.

 We know that the result is a polynomial~$p(x)$ of degree $4n^2$, whose
 coefficient of $x^{2n^2}$ we wish to extract. Hence, the whole computation
 can be done modulo $x^{2n^2+1}$, which does not change the coefficient
 of~$x^{2n^2}$, but which reduces the size of intermediate
 expressions. Alternatively, one can apply the evaluation-interpolation
 technique combined with the Chinese remainder theorem (CRT). Not only
 the degree of~$p(x)$ is known, we also know that it is palindromic, i.e.,
 $p(x)=x^{4n^2}p(1/x)$. Therefore $p(x)$ can be interpolated by using only
 $2n^2+1$ evaluation points. In addition, we can determine an a-priori bound
 on the height of~$p(x)$: let $k$ be the maximal number of nonzero entries in any
 row (or column) of $M$, then the entries of $M^{2n-1}$ are polynomials with
 coefficients at most $k^{2n-2}$ and the height of $p(x)$ is thus bounded
 by $\ell^2k^{2n-2}$.

 Table~\ref{tab:strategies} illustrates the effect of the different strategies on the
 runtime of the computation.
 \begin{table}
 \begin{center}
   \begin{tabular}{c|c|r|r|r}
     \multicolumn{2}{c|}{transfer matrix} &
     mult.\ in $\set Z[x]$ & in $\set Z[x]/(x^{2n^2+1})$ & eval-int.\ + CRT
     \\ \hline \rule{0pt}{12pt}%
     dense  & $M$   & 23.1\,s & 22.3\,s & 69.6\,s \\
     dense  & $M'X$ &  9.1\,s &  8.3\,s & 89.7\,s \\
     sparse & $M$   & 16.3\,s & 15.3\,s &  1.8\,s \\
     sparse & $M'X$ &  1.9\,s &  2.0\,s &  0.6\,s
   \end{tabular}
 \end{center}
 \caption{Effect of different strategies on the runtime, for $n=4$.}\label{tab:strategies}
 \end{table}
 Table~\ref{tab:runtime} shows computation times for $3\leq n\leq 7$. Note that we used
 parallelization for some of the tasks, but the timings are given in CPU time.
 They were measured on Intel Xeon E5-2630v3 processors at 2.4\,GHz.
 As a curiosity, we realized that Mathematica takes about twice as long for computing
 $v_\text{init}^\top\cdot M$ (row vector times matrix) compared to
 computing the equivalent product $M^\top\cdot v_\text{init}$ (transposed
 matrix times column vector).
 \begin{table}
 \begin{center}
   \begin{tabular}{r|r|r|r|r|r}
     \multirow{2}{*}{$n$}
     & \multirow{2}{*}{create states}
     & \multirow{2}{*}{build~$M$}
     & \multicolumn{3}{c}{total time for matrix-vector multiplications}
     \\ \cline{4-6} \rule{0pt}{12pt}%
     & & & mult.\ in $\set Z[x]$ & in $\set Z[x]/(x^{2n^2+1})$ & eval-int.\ + CRT
     \\ \hline \rule{0pt}{12pt}%
     3 & 0.01\,s & 0.39\,s & 0.10\,s & 0.12\,s & 0.06\,s \\
     4 & 0.06\,s & 10\,s  & 1.9\,s & 2.0\,s & 0.6\,s \\
     5 & 0.87\,s & 5\,min & 56\,s & 54\,s & 21\,s \\
     6 & 13.6\,s & 5\,h   & 49\,min & 44\,min & 10\,min \\
     7 & 4\,min & 5\,d   & 29\,h & 25\,h & 5\,h
   \end{tabular}
 \end{center}
 \caption{Runtime observed for various strategies and various problem sizes,
   using the decomposition $M'X$ and the sparse matrix representation.}
 \label{tab:runtime}
 \end{table}
 The timings in the last column of Table~\ref{tab:runtime}
 refer to the minimal number of primes that are
 needed to reconstruct the correct result (which we know already from the two
 previous computations). In order to obtain a provably correct result with the
 CRT approach, one would have to use enough primes to exceed the bound on
 the height of $p(x)$. For example, the result for $n=7$ is approx.\
 $1.335\cdot10^{36}$, therefore requiring $4$ primes of size $2^{31}$, while our
 bound $\ell^2k^{2n-2}$ with $\ell=393878$ and $k=16384$ yields $5.804\cdot10^{61}$,
 corresponding to 7 primes of the same size.

 \section{Further results}\label{sec:4}

 While the gerrymandering problem for a rectangular board of size $m\times n$ is symmetric in $m$ and~$n$,
 the cost of the transfer-matrix method is highly asymmetric.
 As explained above, the cost depends exponentially on the side length that determines the transfer matrix but only
 polynomially on the side length that appears in the exponent.
 Because of this discrepancy, slim rectangular boards are somewhat easier to handle than boards that are
 quadratic or close to quadratic.
 For small values of $m$, it is not too hard to let $n$ grow into the hundreds or even thousands. 

 For fixed~$m$ and varying~$n$, we are interested in the number of solutions to the gerrymandering problem
 for a board of size $m\times n$.
 Obviously there is no solution when both $m$ and $n$ are odd.
 Therefore, for fixed and even $m$, we define $a_n$ as the number of solutions for a board of size $m\times n$,
 and for fixed odd~$m$, we define $a_n$ as the number of solutions for a board of size $m\times 2n$.
 For certain vectors $v_{\text{init}},v_{\text{final}}$ and a certain matrix~$M$ whose entries
 are polynomials in~$x$, we then have $a_n = \frac12[x^{nm/2}](v_{\text{init}}^\top M^{n-1} v_{\text{final}})$
 if $m$ is even and $a_n = \frac12[x^{nm}](v_{\text{init}}^\top M^{2n-1} v_{\text{final}})$ if $m$ is odd.
 We know from linear algebra that the entries of a matrix power $M^n$ are C-finite sequences
 with respect to~$n$, i.e., they satisfy linear recurrences with constant coefficients, or in
 other words, their generating functions are rational.
 An explicit formula is given in Thm.~4.7.2 of Stanley's book~\cite{stanley99}: for a fixed $\ell\times\ell$-matrix $M$
 and any $i,j\in\{1,\dots,\ell\}$, the generating function of the sequence appearing in the $(i,j)$th
 entry of $M^n$ is 
 \[
   (-1)^{i+j} \frac{\det(I_\ell - t\,M)^{[j,i]}}{\det(I_\ell - t\,M)},
 \]
 where the exponent $[j,i]$ indicates the removal of the $j$th row and the $i$th column of
 the matrix $I_\ell-tM$. The generating function for a sequence defined as
 $v_{\text{init}}^\top M^{n-1} v_{\text{final}}$ is just a certain linear combination of such
 rational functions.

 In particular, the rational generating function for the sequence $(v_{\text{init}}^\top M^{n-1} v_{\text{final}})_{n=0}^\infty$
 (or $(v_{\text{init}}^\top M^{2n-1} v_{\text{final}})_{n=0}^\infty$, if $m$ is odd)
 can be explicitly computed from the vectors $v_{\text{init}}, v_{\text{final}}$, and the matrix~$M$.
 At least in principle.
 In practice, for large matrices $M$ involving a symbolic parameter~$x$, computing the determinant
 of $I_\ell-t\,M$, which involves an additional symbolic parameter~$t$, can be a hassle.
 We have managed to compute the rational expressions using evaluation/interpolation techniques
 for $m=3,\dots,7$. The computation is non-rigorous in so far as the number of evaluation points
 was only determined experimentally.
 In Table~\ref{tab:ratsizes}, we summarize their sizes.
 \begin{table}
 \begin{center}
   \begin{tabular}{c|c|c|c|c|c|c|c}
     $m$ & 1 & 2 & 3 & 4 & 5 & 6 & 7 \\\hline
     $\deg_t$ of numerator   & 2 & 4 & 7 & 17 & 36 & 75 & 203 \\
     $\deg_x$ of numerator   & 2 & 4 & 23 & 34 & 190 & 236 & 1425 \\
     monomials in numerator  & 3 & 9 & 76 & 194 & 1955 & 4312 & 55218 \\
     $\deg_t$ of denominator & 2 & 4 & 7 & 17 & 36 & 75 & 203 \\
     $\deg_x$ of denominator  & 2 & 4 & 22 & 34 & 188 & 235 & 1422 \\
     monomials in denominator & 4 & 9 & 76 & 194 & 1935 & 4310 & 55188 \\
   \end{tabular}
 \end{center}
 \caption{Sizes of rational generating functions for various values of~$m$.}\label{tab:ratsizes}
 \end{table}
 The cases $m=1$ and $m=2$ are quite simple, for example, for $m=2$ the generating function is
 \begin{alignat*}1
   &\frac{-t^4 x^4-t^3 x^4-2 t^3 x^3-t^3 x^2+4 t^2 x^2-t x^2+2 t x-t+1}
      {(t-1)^2 (t x^2-1)^2}\\
      &= 1+(x^2+2 x+1)t+(x^4+4 x^3+4 x^2+4 x+1)t^2\\
      &\quad+(x^6+6 x^5+6 x^4+6 x^3+6 x^2+6x+1)t^3+\cdots
 \end{alignat*}
 For $m=3$, the expression is already too big to fit in a line, and as can be seen in the table,
 the sizes increase significantly with respect to~$m$.
 
 The number of colors affects the growth of the expressions even more significantly.
 We have considered a variant of the problem where besides black cells and white cells we also have gray cells.
 The question is then how many ways there are to dissect the $m\times n$ grid into three connected regions,
 with prescribed areas for each color.
 The corresponding rational generating function contains three variables: one marking the length $n$ of the
 board, one marking the area of black cells, and one marking the area of white cells. (There is no need
 for a variable marking the area of the gray cells because we know that the areas of the
 three colors must sum to~$mn$.) 
 We were only able to construct the rational generating function for the case $m=3$.
 Its numerator has degree 29 in $t$ and degree 32 in $x_1$ and $x_2$, it altogether consists of 7939 monomials.
 The denominator has degree 28 in $t$ and degree 30 in $x_1$ and $x_2$, and it consists of 7412 monomials.

 Returning to the case of two colors, it remains to discuss the coefficient extraction operator.
 If $a(x,t)$ is a rational generating function in $x$ and~$t$, viewed as power series in $t$ whose
 coefficients are polynomials in~$x$, extracting the coefficient of $x^{\alpha n}$ (with $\alpha=m/2$
 or $\alpha=m$ depending on whether $m$ is even or odd) from the $n$th
 term of the series is the same as extracting the coefficient of $x^{-1}$ from $x^{-1}a(x,t/x^\alpha)$.
 This can be done by creative telescoping~\cite{Zeilberger91,Chyzak14,Koutschan13a},
 as follows: using computer algebra,
 we can compute polynomials $p_0(t),\dots,p_r(t)$ which only depend on $t$ as well as a rational
 function $b(x,t)$ such that
 \[
   p_0(t)\frac{a(x,t/x^{\alpha})}{x} + \cdots + p_r(t) \frac{d^r}{dt^r}\frac{a(x,t/x^{\alpha})}{x}
   = \frac{d}{dx}b(x,t).
 \]
 Applying $[x^{-1}]$ on both sides, we get zero on the right, because the derivative of a rational
 function cannot have a residue. On the left, observe that taking the residue with respect to $x$
 commutes with the polynomials $p_i(t)$ and the derivations in~$t$. Therefore the series
 $a(t) := \frac12[x^{-1}]\bigl(x^{-1}a(x,t/x^{\alpha})\bigr)$ satisfies the differential equation
 \[
   p_0(t)a(t) + \cdots + p_r(t)\frac{d^r}{dt^r}a(t) = 0.
 \]
 In practice, using the second-named author's implementation~\cite{Koutschan10b}, this approach works nicely for $m=3$ (\oeis{A167242}),
 where we obtain a computer proof of Knuth's result mentioned in the introduction, and for $m=4$ (\oeis{A167247}),
 where we find a differential equation of order~5 and polynomial coefficients of degree~208.
 Using guessing~\cite{Kauers09}, we can also find a differential equation of order~2 with polynomial
 coefficients of degree~96 as well as a polynomial equation of degree~2 with polynomial coefficients
 of degree~51. The correctness of these guessed equations can be proved by showing that the corresponding
 differential operators are right factors of the differential operators obtained via creative telescoping,
 and checking an appropriate number of initial values.
 The equations, as well as the rational generating functions, are
 available at \url{http://www.koutschan.de/data/gerry/}.

 For $m\geq5$, we have not been able to find equations either by creative telescoping or by guessing,
 although Furstenberg's theorem guarantees their existence.
 Following Zeilberger's example, the first-named author (M.K.) will therefore offer a donation
 of \texteuro100 to the OEIS in honor of the person who first manages to find a differential equation
 for some $m\geq5$, either experimentally or rigorously. 

 \section{Acknowledgment}

 We thank Doron Zeilberger for encouraging us to work on this problem and to
 write this paper.  
 We are also indebted to the anonymous referee for checking our manuscript
 very carefully, by even producing an own implementation of our method. This
 process helped to clarify some inaccuracies and thus greatly improved the
 readability of the paper.
 M.K. was supported by the Austrian FWF grant P31571-N32.

	\bigskip
	\hrule
	\bigskip
	
	\noindent 
	2010 \emph{Mathematics Subject Classification}:~Primary 05A15.
        Secondary 68W30, 33F10.
	
	\medskip
	
	\noindent 
	\emph{Keywords}:~transfer-matrix method, polyomino, computer algebra.
	
	\bigskip
	\hrule
	\bigskip
	
	\noindent 
	(Concerned with sequences \oeis{A167242}, \oeis{A167247}, and \oeis{A348456}.)
	
	\bigskip
	\hrule
	\bigskip
	
	\vspace*{+.1in}
	\noindent
	Received ???;
	revised version received ???.
	Published in {\it Journal of Integer Sequences}, ???.
	
	\bigskip
	\hrule
	\bigskip
	
	\noindent
	Return to
	\htmladdnormallink{Journal of Integer Sequences home page}{https://cs.uwaterloo.ca/journals/JIS/}.
	\vskip .1in
 
\end{document}